\documentclass[english,twoside,a4paper,12pt]{amsart}
\usepackage{a4wide}
\usepackage[T1]{fontenc}
\usepackage{amsthm}
\usepackage[mac]{inputenc}
\usepackage{amsfonts,amsmath,epsfig,amscd,amssymb,latexsym}
\theoremstyle{plain}
\newtheorem{theorem}{Theorem}
\newtheorem{lemma}[theorem]{Lemma}

\newtheorem{corollary}[theorem]{Corollary}

\theoremstyle{definition}

\newtheorem{remark}[theorem]{Remark}

\title{Mean square estimate for relatively short exponential sums involving Fourier coefficients of cusp forms}
\author{Anne-Maria Ernvall-Hyt\"onen}\thanks{The author would like to thank prof. Valentin Blomer for fruitful conversations, and for the idea of considering this question. Furthermore, the author would also like to thank Esa Vesalainen for useful conversations during the time working with the problem. The work was funded by the Academy of Finland, grant 138337.}
\address{Department of Mathematics and Statistics, University of Helsinki, Finland}

\begin{document}
\begin{abstract}
We estimate the mean square of a short exponential sum involving Fourier coefficients of a cusp form with a linear twist, a smooth weight function, and a relatively short averaging interval.
\end{abstract}

\maketitle
\section{Introduction}
Let $f(z)=\sum_{n=1}^{\infty}a(n)n^{(\kappa-1)/2}e(nz)$ with $e(x)=e^{2\pi i x}$ be a holomorphic cusp form of weight $\kappa$ with respect to the full modular group. Estimating exponential sums involving Fourier coefficients of cusp forms is a classical question. For the so called long sums, the best possible bound is well-known: Jutila \cite{jutila:ramanujan} has proved that
\[
\sum_{1\leq n\leq M}a(n)e(n\alpha)\ll M^{1/2},
\]
when $\alpha\in \mathbb{R}$. This was an improvement over the classical result by Wilton \cite{wilton}. By the Rankin-Selberg mean value theorem \cite{rankin} this bound is the best possible in the general case, even though for some values of $\alpha$ it is possible to prove considerably better bounds. For instance, for rational values of $\alpha=\frac{h}{k}$, the classical bound is $\ll k^{2/3}M^{1/3+\varepsilon}$

However, the behavior of  short sums, i.e., the sums over an interval $[x,y]$, where $y-x=o(x)$ is much less known. These sums have been investigated for instance by the author and Karppinen in \cite{e&k}. Even though some of the bounds proved there are sharp, it is likely that in many cases, the actual bounds are much smaller than what have been proved.

It is generally a very difficult question to prove good bounds for individual sums. It is much easier to consider the average behavior, namely, to bound expressions of the type
\[
\int_M^{M+\Delta}w(x)\left|\sum_{x\leq n\leq x+y}a(n)e(\alpha n)\right|^2dx.
\] 

Mean squares have been used to deduce the average behavior on various types of objets of interest, including the zeta function and the error term in the divisor problem (for the latter, see Cram\'er \cite{cramer}). For further results on the error term in the divisor problem, see e.g. \cite{heath-browntsang}.

The classical mean square result for long sums can be found in \cite{jutila:lectures}, Theorem 1.2:
\[
\int_1^M\left|\sum_{n\leq x}a(n)e\left(\frac{hn}{k}\right)\right|^2dx=c_2(\kappa)kM^{3/2}+O\left(k^2M^{1+\varepsilon}\right)+O\left(k^{3/2}M^{5/4+\varepsilon}\right)
\]
This result fits very well together with the result we are going to prove in this paper. A corresponding result also holds for sums involving the divisor function twisted with a rational parameter. 

Questions closely connected to the topic of the current paper have also been dealt in \cite{ivic:mean}, where Ivi\'c proves the asymptotic result for $\alpha=0$, $y\ll \sqrt{x}$ and $\Delta=M$, and in \cite{jutila:salomaa}, where Jutila proved an asymptotic result for a mean-square involving a sum of  values of the divisor function with $y\ll x^{1/2}$ and $\Delta\gg M^{1/2}$. In the case $y=\sqrt{x}$ an exponential sum involving Fourier coefficients of a cusp form was dealt in \cite{oma:mean}. There the averaging interval depended on the exponential twist similarly as in the current paper. Very recently, Vesalainen \cite{vesalainen:oskillaatio} proved a mean square result for exponential sums of length at most square root twisted with a rational parameter.
 This leads to the natural question: what happens when the sum is short but longer than of square root length, and the sum has an exponential twist.

We give the following answer to the question:
\begin{theorem}\label{paatulos}
Let $1>\delta> 1/2$, $T\asymp M^{\delta}$, $0<k\ll M^{1/2-\varepsilon'}$ for some fixed positive $\varepsilon'$ and let $h$ and $k$ be co-prime. Let $\Delta=\min(k^2M^{1/2+\varepsilon},M)$, and let $w(x)$ be an infinitely many times differentiable smooth weight function in $\mathbb{R}$ that has support on the interval $[M,M+\Delta]$ satisfying the conditions $w(x)^{(j)}\ll \Delta^{-j}$ for all non-negative $j\in \mathbb{Z}$. Further assume $w(x)=1$ for $M+\Delta'\leq M+\Delta-\Delta'$ for some $\Delta'\asymp \Delta$ and $w(x)=0$ for $x\leq M$ and $x\geq M+\Delta$. Now
\begin{multline*}
\int_{M}^{M+\Delta}w(x)\left|\sum_{x\leq n\leq x+T}a(n)e\left(\frac{h}{k}\right)\right|^2dx=S+O(k^2M^{1+\varepsilon}+\Delta M^{\varepsilon}T^{1/2}k)\\+O\left(\sqrt{|S|(k^2M^{1+\varepsilon}+\Delta M^{\varepsilon}T^{1/2}k)}\right),
\end{multline*}
where
\begin{align*}
S=\frac{2k}{\pi^2}\sum_{n\leq\min(M^2T^{-2},M^{2}T^{-1}\Delta^{-1})}|a(n)|^2n^{-3/2}\int_M^{M+\Delta}w(x)x^{1/2}\sin^2\left(2\pi \frac{\sqrt{n(x+T)}}{k}-2\pi\frac{\sqrt{nx}}{k}\right)\\ \times\sin^2\left(2\pi\frac{\sqrt{n(x+T)}}{k}+2\pi\frac{\sqrt{nx}}{k}-\frac{\pi}{4}\right)dx.
\end{align*}
Furthermore,
\[
S\ll \begin{cases}\Delta TM^{\varepsilon}&\textrm{when }k\gg TM^{-1/2} \\ \Delta M^{1/2+\varepsilon}k &\textrm{otherwise.}\end{cases}
\]
\end{theorem}

\begin{remark}
We have to assume that $1\leq k \ll M^{1/2-\varepsilon'}$ for some positive $\varepsilon'$. If we let $k$ to be larger, then the error term in the truncated Voronoi type summation formula would be too large, i. e., we would not get a better average bound than the one obtained from Jutila's estimate for long sums \cite{jutila:ramanujan} using the triangle inequality.
\end{remark}

Notice that if we are only interested in the upper bound, we can easily deduce a result without the weight function, since by choosing the weight function $w^{\star}(x)$ to be a smooth weight function supported on the interval $\left[M-\frac{\Delta}{2},M+\frac{3\Delta}{2}\right]$, and obtaining value $w^{\star}(x)=1$ on $\left[M,\Delta\right]$, we have
\[
\int_M^{M+\Delta}\left|\sum_{x\leq n\leq x+T}a(n)e\left(\frac{h}{k}\right)\right|^2dx\leq \int_{M-\Delta/2}^{M+3\Delta/2}w^{\star}(x)\left|\sum_{x\leq n\leq T}a(n)e\left(\frac{h}{k}\right)\right|^2dx,
\]
and hence, noticing that the term $\Delta M^{\varepsilon} T^{1/2}k$ is always smaller than the other terms, we obtain the following corollary:
\begin{corollary}
Let $1>\delta> 1/2$, $T\asymp M^{\delta}$, $0<k\ll M^{1/2-\varepsilon'}$ for some fixed positive $\varepsilon'$ and let $h$ and $k$ be co-prime. Let $\Delta=\min(k^2M^{1/2+\varepsilon},M)$. Now
\[
\int_{M}^{M+\Delta}\left|\sum_{x\leq n\leq x+T}a(n)e\left(\frac{h}{k}\right)\right|^2dx\ll \min\left(\Delta TM^{\varepsilon},\Delta M^{1/2+\varepsilon}k\right)+O(k^2M^{1+\varepsilon}).
\]
\end{corollary}

Writing this corollary according to various values of $k$‚ we obtain the following bounds:
\begin{corollary}
Let $1>\delta> 1/2$, $T\asymp M^{\delta}$, $0<k\ll M^{1/2-\varepsilon'}$ for some fixed positive $\varepsilon'$ and let $h$ and $k$ be co-prime. Let $\Delta=\min(k^2M^{1/2+\varepsilon},M)$. Now
\[
\int_{M}^{M+\Delta}\left|\sum_{x\leq n\leq x+T}a(n)e\left(\frac{h}{k}\right)\right|^2dx\ll \begin{cases}\Delta M^{1/2+\varepsilon}k& \textrm{if $k\ll TM^{-1/2}$}\\ \Delta TM^{\varepsilon} & \textrm{if $k\gg TM^{-1/2}$ and $k\ll T^{1/2}$}\\ \Delta M^{\varepsilon}k^2& \textrm{if $k\gg TM^{-1/2}$ and $k\gg T^{1/2}$}.\end{cases}
\]
\end{corollary}

The main advantage in this theorem is the relatively short averaging interval compared to the length of the sum, when the value of $k$ is small. The averaging interval is actually similar to the one in \cite{oma:mean}. However, in the current paper, the length of the sum can be much longer than the averaging interval unlike there. In particular, we obtain the following corollary:
\begin{corollary}
We have
\[
\int_M^{M+\Delta}\left|\sum_{x\leq n\leq x+T}a(n)\right|^2dx\ll M^{1/2+\varepsilon}\Delta,
\]
where $\Delta \gg M^{1/2+\varepsilon'}$, and where $\varepsilon>0$ and $\varepsilon'>0$ can be chosen to be arbitrarily small fixed positive numbers.
\end{corollary}

This corollary shows that already on a short averaging interval the average behavior of a nearly long sum (i. e., a sum of length $x^{1-\varepsilon}$) has the conjectured size $x^{1/4+\varepsilon}$ in average.

In order to prove Theorem \ref{paatulos}, we need some lemmas, which will be presented in the following section. In the final section, we will have the proof of Theorem \ref{paatulos}.

In the following, the $\varepsilon$'s will be positive, not necessarily equal. The constants implied by symbols $\ll, \gg$ and $O()$ do not depend on $M$ or $k$, but they do depend on $\varepsilon$, $\delta$, on the properties of the weight function, etc.

\section{Lemmas}

The first lemma is partial integration (see e.g. \cite{jutimoto:acta} Lemma 6):
\begin{lemma}\label{jutilamotohashi}
Let $A$ be a $P\geq 0$ times differentiable function which is
compactly supported in a finite interval $[a,b]$. Assume also that
there exist two quantities $A_0$ and $A_1$ such that for any
non-negative integer $\nu\leq P$ and for any $x\in [a,b]$,
\[
A^{(\nu)}(x)\ll A_0A_1^{-\nu}.
\]
Moreover, let $B$ be a function which is real-valued on $[a,b]$,
and regular throughout the complex domain composed of all points
within the distance $\varrho$ from the interval; and assume that
there exists a quantity $B_1$ such that
\[
0<B_1\ll \left|B'(x)\right|
\]
for any point $x$ in the domain. Then we have
\[
\int_{a}^{b}A(x)e\left(B(x)\right)d x \ll
A_0\left(A_1B_1\right)^{-P}\left(1+\frac{A_1}{\varrho}\right)^P\left(b-a\right).
\]
\end{lemma}

The following lemmas can be proved using partial integration or the previous lemma. The details are similar to those in \cite{oma:mean}.
\begin{lemma}\label{erifunktiot} Let $1\leq m,n\leq M$. Then
\[
\int_M^{M+\Delta} w(x)x^{1/2}e\left(\pm \left(2\frac{\sqrt{nT_1(x)}}{k} +2\frac{\sqrt{mT_2(x)}}{k}\right)\right)d x\ll k^P\left(\sqrt{n}+\sqrt{m}\right)^{-P}\Delta^{1-P}M^{P/2+1/2},
\]
where $T_1(x)$ and $T_2(x)$ are $x$ or $x+T$ (not necessarily but possibly the same).
\end{lemma}
\begin{lemma}\label{erimerkki} Let $1\leq m,n\leq M$.Then
\[
\int_M^{M+\Delta} w(x)x^{1/2}e\left(\pm 2\frac{\sqrt{nT(x)}}{k}\mp 2\frac{\sqrt{mT(x)}}{k}\right)d x\ll k^P\left|\sqrt{n}-\sqrt{m}\right|^{-P}\Delta^{1-P}M^{P/2+1/2},
\]
where $T(x)=x$ or $T(x)=x+T$.
\end{lemma}
\begin{lemma}\label{hankalin} Let $1\leq m,n\leq M$.Then
\begin{multline*}
\int_M^{M+\Delta}x^{1/2}w(x)e\left(\pm2\frac{\sqrt{m(x+T)}}{k}\mp2\frac{\sqrt{nx}}{k}\right)d x\\ \ll \max_{x\in [M,M+\Delta]}k^P\Delta^{1-P}\left|\frac{\sqrt{m}}{\sqrt{1+Tx^{-1}}}-\sqrt{n}\right|^{-P}M^{P/2+1/2}.
\end{multline*}
\end{lemma}
When $\frac{\sqrt{m}}{\sqrt{1+Tx^{-1}}}-\sqrt{n}=0$ for some $x\in [M,M+\Delta]$, Lemma \ref{hankalin} does not give any information. 

\section{Proof of Theorem \ref{paatulos}}
Using a Voronoi-type summation formula (see \cite{jutila:lectures} Theorem 1.2 with $x\asymp M$ and the choice $N=M$, we have
\[
\sum_{n\leq x}a(n)e\left(\frac{h}{k}n\right)=\left(\pi \sqrt{2}\right)^{-1}x^{1/4}k^{1/2}\sum_{n\leq M}a(n)e\left(-\frac{\bar{h}}{k}n\right)n^{-3/4}\cos \left(4\pi\frac{\sqrt{nx}}{k}-\frac{\pi}{4}\right)+O\left(M^{\varepsilon}k\right)
\]
Squaring and integration over the error term gives a total contribution of at most
\[
\Delta M^{\varepsilon}k^2\ll M^{1+\varepsilon}k^2.
\]
We may thus forget it for awhile because the contribution from the cross-terms between the error term and the other terms can be taken into account using the Cauchy-Schwarz inequality. Denote
\[
S(x,M_1,M_2)=\left(\pi \sqrt{2}\right)^{-1}x^{1/4}k^{1/2}\sum_{M_1<n\leq M_2}a(n)e\left(-\frac{\bar{h}}{k}n\right)n^{-3/4}\cos \left(4\pi\frac{\sqrt{nx}}{k}-\frac{\pi}{4}\right).
\]
Hence, the expression we need to consider is
\[
\int_M^{M+\Delta}w(x)\left|S(x+T,0,M)-S(x,0,M)\right|^2dx.
\]
Let us first split the summation at $M'$:
\begin{align*}
\int_M^{M+\Delta}&w(x)\left|S(x+T,0,M)-S(x,0,M)\right|^2dx\\&=\int_M^{M+\Delta}w(x)\left|S(x+T,0,M')-S(x,0,M')\right|^2dx\\&+\int_M^{M+\Delta}w(x)\left|S(x+T,M',M)-S(x,M',M)\right|^2\\ & +O\left(\int_M^{M+\Delta}w(x)\left|S(x+T,M',M)-S(x,M',M)\right|\left|S(x+T,0,M')-S(x,0,M')\right|dx\right),\end{align*}
and choose $M'$ in the following way:
\[
M'=\begin{cases}M^2T^{-2}& \textrm{if } T\gg \Delta\\ M^2T^{-1}\Delta^{-1}& \textrm{if }\Delta \gg T.\end{cases}
\]
Notice that if $\Delta=M$, then $\Delta\gg T$, and if $\Delta=M^{1/2+\varepsilon}k^2$, then the condition $\Delta \gg T$ is equivalent to $M^{1/2+\varepsilon}k^2\gg T$, so $k\gg T^{1/2}M^{-1/4-\varepsilon}$.
Let us concentrate on treating the first two terms, and then use the Cauchy-Schwarz inequality to estimate the last term. Let us start with the second term. We have
\begin{multline*}
\int_M^{M+\Delta}w(x)\left|S(x+T,M',M)-S(x,M',M)\right|^2dx\\\ll \int_M^{M+\Delta}w(x)\left|S(x+T,M',M)\right|^2dx+\int_M^{M+\Delta}w(x)\left|S(x,M',M)\right|^2dx\end{multline*}
Let us treat the second term as the first one can be treated similarly. We have
\begin{multline*}
\left|S(x,M',M)\right|^2\\=2\pi^2x^{1/2}k\sum_{M'<m,n\leq M}a(n)\overline{a(m)}e\left(\frac{\bar{h}}{k}(m-n)\right)(mn)^{-3/4}\cos \left(4\pi\frac{\sqrt{nx}}{k}-\frac{\pi}{4}\right)\cos \left(4\pi\frac{\sqrt{mx}}{k}-\frac{\pi}{4}\right)\end{multline*}
The terms with $m=n$ yield a total contribution
\[
2\pi^2k\sum_{M'<n\leq M}|a(n)|^2n^{-3/2}\int_M^{M+\Delta}w(x)x^{1/2}\cos^2\left(4\pi\frac{\sqrt{nx}}{k}-\frac{\pi}{4}\right)dx\ll kM^{1/2}\Delta (M')^{-1/2+\varepsilon}.
\]
When $T\gg \Delta$, we have $M'=M^2T^{-2}$ and hence
\[
kM^{1/2}\Delta (M')^{-1/2+\varepsilon}\ll kM^{1/2}\Delta M^{-1+\varepsilon}T^1=k\Delta TM^{-1/2+\varepsilon}\ll k\Delta M^{\varepsilon}T^{1/2}.
\]
When $T\ll \Delta$, we have $M'=M^2T^{-1}\Delta^{-1}$, and hence
\[
kM^{1/2}\Delta (M')^{-1/2+\varepsilon}\ll kM^{1/2}\Delta M^{-1+\varepsilon}T^{1/2}\Delta^{1/2}\ll k\Delta M^{\varepsilon}T^{1/2}.
\]
When $m\ne n$, we need to estimate integrals
\[
\int_M^{M+\Delta}w(x)x^{1/2}e\left(\pm 2(\sqrt{n}-\sqrt{m})\frac{\sqrt{x}}{k}\right)dx
\]
and
\[
\int_M^{M+\Delta}w(x)x^{1/2}e\left(\pm 2(\sqrt{n}+\sqrt{m})\frac{\sqrt{x}}{k}\right)dx.
\]
We only consider the first one as the second one is similar but simpler. We use Lemma \ref{erimerkki} to bound the integral. We obtain
\[
\int_M^{M+\Delta}w(x)x^{1/2}e\left(\pm 2(\sqrt{n}-\sqrt{m})\frac{\sqrt{x}}{k}\right)dx\ll \Delta^{1-P}k^{P}M^{1/2+P/2}|\sqrt{n}-\sqrt{m}|^{-P}.
\]
The sum over these estimates can be made as small as desired if
\[
\Delta^{-1}(\sqrt{n}-\sqrt{m})^{-1}k\sqrt{M}\ll M^{-\varepsilon}.
\]
This is the case for all except $\asymp \frac{kM^{1/2+\varepsilon}}{\Delta}\sqrt{n}$ values around each $n$ (as long as $n$ is sufficiently large such that there are values $m\ne n$ in this neighborhood, which is the case when $n\gg \frac{\Delta^2}{k^2M^{1+\varepsilon}}$). Let us use absolute values in these cases to bound the integral, and sum over these values:
\begin{multline*}
2\pi^2k\sum_{\substack{\frac{\Delta^2}{k^{2}M^{1+\varepsilon}}<n\leq M\\|n-m|\ll \frac{kM^{1/2+\varepsilon}}{\Delta}\sqrt{n}}}a(n)\overline{a(m)}(nm)^{-3/4}e\left((m-n)\frac{\bar{h}}{k}\right)\int_M^{M+\Delta}w(x)x^{1/2}e\left(\pm 2(\sqrt{n}-\sqrt{m})\frac{\sqrt{x}}{k}\right)dx\\
\ll k\sum_{\substack{\frac{\Delta^2}{k^{2}M^{1+\varepsilon}}<n\leq M\\|n-m|\ll \frac{kM^{1/2+\varepsilon}}{\Delta}\sqrt{n}}}(nm)^{\varepsilon-3/4}\Delta M^{1/2} \ll kM^{1/2}\Delta \sum_{n}n^{-3/2+\varepsilon} \frac{kM^{1/2+\varepsilon}}{\Delta}\sqrt{n}\ll k^2M^{1+\varepsilon}.\end{multline*}

We have now derived the total contribution from the terms with $m,n\gg M'$ and the error term in the Voronoi summation formula to be $\ll k^2M^{1+\varepsilon}+k\Delta M^{\varepsilon}T^{1/2}$.

Let us now move to considering the part
\[
\int_M^{M+\Delta}w(x)\left|S(x+T,0,M')-S(x,0,M')\right|^2dx.
\]
This part is technically somewhat more challenging because it is not sufficient to just use the triangle inequality, but we need to get some cancellation on the diagonal terms. Therefore, we will also meet some integrals that are more difficult to consider than the integrals treated in the first part of the proof.

We have
\begin{align*}
S(x+T,0,M')& =\frac{k^{1/2}}{\pi\sqrt{2}}\left(x+T\right)^{1/4}\sum_{n\leq M'}a(n)e\left(-\frac{\bar{h}}{k}n\right)n^{-3/4}\cos \left(4\pi\frac{\sqrt{n(x+T)}}{k}-\frac{\pi}{4}\right)\\
&=\frac{k^{1/2}}{\pi\sqrt{2}}x^{1/4}\sum_{n\leq M'}a(n)e\left(-\frac{\bar{h}}{k}n\right)n^{-3/4}\cos \left(4\pi\frac{\sqrt{n(x+T)}}{k}-\frac{\pi}{4}\right)\\ &-\frac{k^{1/2}}{\pi\sqrt{2}}\left(\left(x+T\right)^{1/4}-x^{1/4}\right)\sum_{n\leq M'}a(n)e\left(-\frac{\bar{h}}{k}n\right)n^{-3/4}\cos \left(4\pi\frac{\sqrt{n(x+T)}}{k}-\frac{\pi}{4}\right)\\&=\frac{k^{1/2}}{\pi\sqrt{2}}x^{1/4}\sum_{n\leq M'}a(n)e\left(-\frac{\bar{h}}{k}n\right)n^{-3/4}\cos \left(4\pi\frac{\sqrt{n(x+T)}}{k}-\frac{\pi}{4}\right)\\ & +O\left(k^{1/2}(M')^{1/4+\varepsilon}M^{-3/4}T\right).\end{align*}
Thus, we have
\begin{multline*}
S(x+T,0,M')-S(x,0,M')\\=\frac{k^{1/2}}{\pi\sqrt{2}}x^{1/4}\sum_{n\leq M'}a(n)e\left(-\frac{\bar{h}}{k}n\right)n^{-3/4}\left(\cos \left(4\pi\frac{\sqrt{n(x+T)}}{k}-\frac{\pi}{4}\right)-\cos \left(4\pi\frac{\sqrt{nx}}{k}-\frac{\pi}{4}\right)\right)\\+O\left(k^{1/2}(M')^{1/4+\varepsilon}M^{-3/4}T\right).\end{multline*}
Let us first treat the error term and then concentrate on the main term. Squaring and integrating over the main term gives the total contribution
\[
\Delta k(M')^{1/2+\varepsilon}M^{-3/2}T^2.
\]
If $T\gg \Delta$, we have $M'=M^2T^{-2}$, and hence
\[
\Delta k(M')^{1/2+\varepsilon}M^{-3/2}T^2\ll \Delta kM^{1+\varepsilon}T^{-1}M^{-3/2}T^2=kM^{-1/2+\varepsilon}\Delta T\ll k\Delta M^{\varepsilon}T^{1/2}.
\]
If $\Delta \gg T$, we have $M'=M^2T^{-1}\Delta^{-1}$, and hence
\[
\Delta k(M')^{1/2+\varepsilon}M^{-3/2}T^2\ll \Delta kM^{1+\varepsilon}T^{-1/2}\Delta^{-1/2}M^{-3/2}T^2=kM^{-1/2+\varepsilon}\Delta^{1/2} T^{3/2}\ll k\Delta M^{\varepsilon}T^{1/2}.
\]
In both cases, the contribution coming from the error term is at most $k\Delta M^{\varepsilon}T^{1/2}$. Let us now move to the main term.
\begin{multline*}
\int_M^{M+\Delta}w(x)\left|\frac{k^{1/2}}{\pi\sqrt{2}}x^{1/4}\sum_{n\leq M'}a(n)e\left(-\frac{\bar{h}}{k}n\right)n^{-3/4}\right.\\\times \left.\left(\cos \left(4\pi\frac{\sqrt{n(x+T)}}{k}-\frac{\pi}{4}\right)-\cos \left(4\pi\frac{\sqrt{nx}}{k}-\frac{\pi}{4}\right)\right)\right|^2dx\\=\frac{k}{2\pi^2}\sum_{m,n\leq M'}a(n)\overline{a(m)}e\left((m-n)\frac{\bar{h}}{k}\right)(nm)^{-3/4}\int_M^{M+\Delta}w(x)x^{1/2}\\ \times \left(\cos \left(4\pi\frac{\sqrt{n(x+T)}}{k}-\frac{\pi}{4}\right)-\cos \left(4\pi\frac{\sqrt{nx}}{k}-\frac{\pi}{4}\right)\right)\\ \times \left(\cos \left(4\pi\frac{\sqrt{m(x+T)}}{k}-\frac{\pi}{4}\right)-\cos \left(4\pi\frac{\sqrt{mx}}{k}-\frac{\pi}{4}\right)\right)dx.
\end{multline*}
Let us first look at the case with $n=m$. Now
\begin{multline*}
\left(\cos \left(4\pi\frac{\sqrt{n(x+T)}}{k}-\frac{\pi}{4}\right)-\cos \left(4\pi\frac{\sqrt{nx}}{k}-\frac{\pi}{4}\right)\right)^2\\=4\sin^2\left(2\pi \frac{\sqrt{n(x+T)}}{k}-2\pi\frac{\sqrt{nx}}{k}\right)\sin^2\left(2\pi\frac{\sqrt{n(x+T)}}{k}+2\pi\frac{\sqrt{nx}}{k}-\frac{\pi}{4}\right)
\end{multline*}
When $n\ll \frac{xk^2}{T^2}$, we have
\[
\frac{\sqrt{n(x+T)}}{k}-\frac{\sqrt{nx}}{k}\ll \frac{\sqrt{n}T}{k\sqrt{x}}\ll 1.
\]
Hence, if $\frac{Mk^2}{T^2}\gg 1$, and since by our choice of $M'$, we always have $M'\gg \frac{Mk^2}{T^2}$, we can split the diagonal sum into two. The first part gives the bound
\begin{multline*}
\frac{2k}{\pi^2}\sum_{n\ll \frac{Mk^2}{T^2}}|a(n)|^2n^{-3/2}\int_M^{M+\Delta}w(x)x^{1/2}\sin^2\left(2\pi \frac{\sqrt{n(x+T)}}{k}-2\pi\frac{\sqrt{nx}}{k}\right)\\ \times\sin^2\left(2\pi\frac{\sqrt{n(x+T)}}{k}+2\pi\frac{\sqrt{nx}}{k}-\frac{\pi}{4}\right)dx\\ \ll k\sum_{n\ll \frac{Mk^2}{T^2}}|a(n)|^2n^{-3/2}M^{1/2}\left(\frac{\sqrt{n}T}{k\sqrt{M}}\right)^2\Delta \ll T^2\Delta k^{-1}M^{-1/2}\sum_{n\ll \frac{Mk^2}{T^2}}n^{\varepsilon-1/2}\\ \ll \Delta TM^{\varepsilon}.
\end{multline*}
The second part gives the bound
\begin{multline*}
\frac{2k}{\pi^2}\sum_{n\gg \frac{Mk^2}{T^2}}|a(n)|^2n^{-3/2}\int_M^{M+\Delta}w(x)x^{1/2}\sin^2\left(2\pi \frac{\sqrt{n(x+T)}}{k}-2\pi\frac{\sqrt{nx}}{k}\right)\\ \times\sin^2\left(2\pi\frac{\sqrt{n(x+T)}}{k}+2\pi\frac{\sqrt{nx}}{k}-\frac{\pi}{4}\right)dx\\ \ll k\sum_{n\gg \frac{Mk^2}{T^2}}|a(n)|^2n^{-3/2}M^{1/2}\Delta \ll \Delta M^{1/2}k\left(\frac{Mk^2}{T^2}\right)^{-1/2+\varepsilon}\ll \Delta M^{\varepsilon}T.
\end{multline*}
On the other hand, if $\frac{Mk^2}{T^2}\ll 1$, we have
\begin{multline*}
\frac{2k}{\pi^2}\sum_{n\geq 1}|a(n)|^2n^{-3/2}\int_M^{M+\Delta}w(x)x^{1/2}\sin^2\left(2\pi \frac{\sqrt{n(x+T)}}{k}-2\pi\frac{\sqrt{nx}}{k}\right)\\ \times\sin^2\left(2\pi\frac{\sqrt{n(x+T)}}{k}+2\pi\frac{\sqrt{nx}}{k}-\frac{\pi}{4}\right)dx\\ \ll k\sum_{n\geq 1}|a(n)|^2n^{-3/2}M^{1/2}\Delta \ll \Delta M^{1/2}k.
\end{multline*}
We can now move to the terms with $n\ne m$. We split the cosine product into exponential terms, and thus, we have to treat integrals of the type
\[
\int_M^{M+\Delta}w(x)x^{1/2}e\left(\pm 2\left(\frac{\sqrt{nT_1(x)}}{k}\pm \frac{\sqrt{nT_2(x)}}{k}\right)\right),
\]
where $T_1(x)$ and $T_2(x)$ are $x$ or $x+T$, possibly but not necessarily the same. These integrals have been treated in Lemmas \ref{erifunktiot}-\ref{hankalin}. Lemma \ref{erifunktiot} always gives a good enough bound. Lemma \ref{erimerkki} gives a good enough bound, namely, the integral can be made as small as desired, unless $|\sqrt{n}-\sqrt{m}|\ll M^{1/2+\varepsilon}k\Delta^{-1}$. As in the first part of the proof, there are $\asymp \frac{kM^{1/2+\varepsilon}}{\Delta}\sqrt{n}$ values of $m$ around each $n$ for which $|\sqrt{n}-\sqrt{m}|\ll M^{1/2+\varepsilon}k\Delta^{-1}$. The contribution coming from these terms is again
\[
\ll k^2M^{1+\varepsilon}.
\]
We may now move to the cases when using Lemma \ref{hankalin} fails. Notice first that writing
\[
f_{n,m}(x)=2\frac{\sqrt{n(x+T)}}{k}-2\frac{\sqrt{mx}}{k},
\]
we have
\[
f_{n_0,m_0}'(x)=\frac{\sqrt{n_0}}{\sqrt{x+T}k}-\frac{\sqrt{m_0}}{\sqrt{x}k}=0,
\]
when $n_0=m_0\frac{y+T}{y}$ for some $y\in [M,M+\Delta]$. Next we want to show that when $n$ lies outside the interval $I_{m,c}=\left[m\left(1+\frac{T}{M+\Delta}\right)-c,m\left(1+\frac{T}{M}\right)+c\right]$ for $c\asymp \frac{M^{1/2+\varepsilon}k}{\Delta}\sqrt{m}$, then
\[
f_{n,m}'\gg M^{\varepsilon'}\Delta^{-1}
\]
for some $\varepsilon'$, which is just enough that we can use partial integration Lemma \ref{jutilamotohashi} (and therefore, Lemma \ref{hankalin}) to estimate the integral. Write now
\[
n=m\left(1+\frac{T}{M+\Delta}\right)-c.
\]
Then
\begin{align*}
\left|f_{n,m}'(x)\right|&=\left|\frac{\sqrt{m\left(1+\frac{T}{M+\Delta}\right)-c}}{\sqrt{x+T}k}-\frac{\sqrt{m}}{\sqrt{x}k}\right|=\left|\frac{\sqrt{m}}{\sqrt{x}k}\left(\frac{\sqrt{1+\frac{T}{M+\Delta}-\frac{c}{m}}}{\sqrt{1+\frac{T}{x}}}-1\right)\right|\\ & \asymp \left|\frac{\sqrt{m}}{\sqrt{x}k}\left(\frac{T}{M+\Delta}-\frac{T}{x}-\frac{c}{m}\right)\right|\gg \frac{c}{\sqrt{xm}k}\asymp M^{\varepsilon}\Delta^{-1}.
\end{align*}
We can show a similar estimate for $n=m\left(1+\frac{T}{M}\right)+c$. Notice that when $m\ll \frac{\Delta^2}{M^{1+\varepsilon}k^2}$, this interval contains at most a constant number of terms. Furthermore, by our choice of $M'$, there is at most one $n$ for each $m$ satisfying the condition $n=m\left(1+\frac{T}{y}\right)$ for some $y\in [M,M+\Delta]$, since the length of the interval $I_{m}=\left[m\left(1+\frac{T}{M+\Delta}\right),m\left(1+\frac{T}{M}\right)\right]$ is
\[
m\left(1+\frac{T}{M}\right)-m\left(1+\frac{T}{M+\Delta}\right)\asymp\frac{mT\Delta}{M^2}.
\]
All the integers on the interval $\left[m\left(1+\frac{T}{M+\Delta}\right),m\left(1+\frac{T}{M}\right)\right]$ are greater than $m$. When the distance $mTM^{-1}$ between the largest number on the interval and $m$ is at most $o(1)$, there cannot be any numbers on the interval. This is the case when $m=o(MT^{-1})$.

In particular, when both $m=o\left( \frac{\Delta^2}{M^{1+\varepsilon}k^2}\right)$ and $m=o(MT^{-1})$, the sum is empty. Therefore, it suffices to estimate the sums $m\gg \frac{\Delta^2}{M^{1+\varepsilon}k^2}$ and $m\gg MT^{-1}$. We will just take absolute values of the integrals. The contribution coming from interval $I_{c,m}\setminus I_m$ is
\begin{align*}\ll k\sum_{\substack{m\gg \frac{\Delta^2}{M^{1+\varepsilon}k^2}\\n\in I_{m,c}\setminus I_m}}|a(n)||a(m)|(nm)^{-3/4}\Delta M^{1/2}\ll k^2M^{1+\varepsilon}.
\end{align*}
The contribution coming from interval $I_m$ is
\[
\ll k\sum_{\substack{m\gg MT^{-1}\\n\in I_m}}|a(n)||a(m)|(nm)^{-3/4}\Delta M^{1/2}\ll k\Delta M^{1/2}\sum_{m\gg MT^{-1}} m^{\varepsilon-3/2}\ll k\Delta M^{\varepsilon}T^{1/2}
\]
We have now derived the contribution from the terms $m,n\ll M'$ to be the contribution from the main term and an error of size at most $\ll k^2M^{1+\varepsilon}+k\Delta M^{\varepsilon}T^{1/2}$. Using the Cauchy-Schwarz inequality completes the proof.


\begin{thebibliography}{10}

\bibitem{cramer}
H.~Cram\'er.
\newblock \"uber zwei {S}\"atze von {H}errn {G}. {H}. {H}ardy.
\newblock {\em Math. Z.}, 15:201--210, 1922.

\bibitem{oma:mean}
Anne-Maria Ernvall-Hyt{\"o}nen.
\newblock On the mean square of short exponential sums related to cusp forms.
\newblock {\em Funct. Approx. Comment. Math.}, 45(part 1):97--104, 2011.

\bibitem{e&k}
Anne-Maria Ernvall-Hyt{\"o}nen and Kimmo Karppinen.
\newblock On short exponential sums involving {F}ourier coefficients of
  holomorphic cusp forms.
\newblock {\em Int. Math. Res. Not. IMRN}, (10):Art. ID. rnn022, 44, 2008.

\bibitem{heath-browntsang}
D.~R. Heath-Brown and K.~Tsang.
\newblock Sign changes of ${E}(t)$, ${\Delta}(x)$ and ${P}(x)$.
\newblock {\em J. of Number Theory}, 49:73--83, 1994.

\bibitem{ivic:mean}
Aleksandar Ivi{\'c}.
\newblock On the divisor function and the {R}iemann zeta-function in short
  intervals.
\newblock {\em Ramanujan J.}, 19(2):207--224, 2009.

\bibitem{jutila:lectures}
M.~Jutila.
\newblock {\em Lectures on a {M}ethod in the {T}heory of {E}xponential {S}ums},
  volume~80 of {\em Tata Institute of Fundamental Research Lectures on
  Mathematics and Physics}.
\newblock Published for the Tata Institute of Fundamental Research, Bombay,
  1987.

\bibitem{jutila:ramanujan}
M.~Jutila.
\newblock On exponential sums involving the {R}amanujan function.
\newblock {\em Proc. Indian Acad. Sci. Math. Sci.}, 97(1-3):157--166 (1988),
  1987.

\bibitem{jutila:salomaa}
Matti Jutila.
\newblock On the divisor problem for short intervals.
\newblock {\em Ann. Univ. Turku. Ser. A I}, (186):23--30, 1984.
\newblock Studies in honour of Arto Kustaa Salomaa on the occasion of his
  fiftieth birthday.

\bibitem{jutimoto:acta}
Matti Jutila and Yoichi Motohashi.
\newblock Uniform bound for {H}ecke {$L$}-functions.
\newblock {\em Acta Math.}, 195:61--115, 2005.

\bibitem{rankin}
R.~A. Rankin.
\newblock Contributions to the theory of {R}amanujan's function $\tau(n)$ and
  similar arithmetical functions ii. {T}he order of {F}ourier coefficients of
  integral modular forms.
\newblock {\em Proc. Cambridge Philos. Soc.}, 35:357--372, 1939.

\bibitem{vesalainen:oskillaatio}
E.V. Vesalainen.
\newblock Moments and oscillations of exponential sums related to cusp forms.
\newblock Submitted.

\bibitem{wilton}
J.~R. Wilton.
\newblock A note on {R}amanujan's arithmetical function $\tau(n)$.
\newblock {\em Proc. Cambridge Philos. Soc.}, 25(II):121--129, 1929.

\end{thebibliography}
\end{document}